\documentclass[preprints,article,accept,moreauthors,pdftex]{Definitions/mdpi} 

\usepackage{bm}
\usepackage{bbm}
\usepackage{listings}
\usepackage{subfig}
\usepackage{amsmath}
\usepackage{amssymb}
\usepackage{booktabs}
\usepackage{hyperref}
\usepackage{fourier}
\usepackage[T1]{fontenc}
\usepackage[]{avant}
\usepackage{longtable}
\usepackage[ruled, lined, linesnumbered, commentsnumbered, longend]{algorithm2e}

\hypersetup{colorlinks,citecolor=blue}

\firstpage{1} 
\pubvolume{xx}
\issuenum{1}
\articlenumber{1}
\pubyear{2023}
\copyrightyear{2023}
\history{}

\newcommand*\diff{\mathop{}\!\kern0pt\mathrm{d}}

\Title{Instabilities of explicit finite difference schemes with ghost points on the diffusion equation}
\Author{Fabien Le Floc'h}

\AuthorNames{Fabien Le Floc'h}

\address{fabien@2ipi.com}
\abstract{Ghost, or fictitious points allow to capture boundary conditions that are not located on the finite difference grid discretization. We explore in this paper the impact of ghost points on the stability of the explicit Euler finite difference scheme in the context of the diffusion equation. In particular, we consider the case of a one-touch option under the Black-Scholes model. The observations and results are however valid for a much wider range of financial contracts and models.}

\begin{document}
	
	\section{Introduction}
Under the Black-Scholes model, the price of barrier options financial contracts is strongly impacted by the term-structure of interest rates, dividends and volatilities.  In order to take those into account, the standard practice is to rely on a numerical scheme to price those, as the closed form formulae for the Black-Scholes model assume constant rates, dividend yield and volatility. Finite difference methods are an efficient technique on this kind of problem, since the dimensionality of the partial differential equation involved is low. This stays true under the local volatility model or under the Heston stochastic volatility model.

To make the problem concrete, we will consider in this paper the example of the one-touch option, which pays 1 USD if the underlying asset spot price $S$ moves over the specified barrier level $L^+$ at some time before the maturity date $T$ and 0 USD otherwise. The observations and results stay relevant in more general settings, for example for knock-in barrier options, or more exotic derivative products. The Black-Scholes PDE for the contract value $V(S,t)$ reads \citep{Shreve04}
\begin{equation}
\frac{\partial V}{\partial t} + \frac{1}{2}\sigma^2(t)S^2 \frac{\partial^2 V}{\partial S^2} + (r(t)-q(t))\frac{\partial V}{\partial S} - r(t)V = 0
\end{equation}
where $r$ is the interest rate, $q$ the dividend yield and $\sigma$ the Black-Scholes volatility.

At $S=0$, we let the price be linear. The one-touch contract imposes a Dirichlet boundary condition $V(L^+,t)=1$ for $t \leq T$ and the initial condition reads $V(S,T)=1_{S \geq L^+}$. 

 The most popular way to enforce the boundary condition at $L^+$ is to make sure that there is a grid point at $L^+$. For a  one-touch option, we only need to set the grid upper bound to $L^+$. More exotic structures may involve a set of barrier levels active at different dates. It is then less obvious to properly define the grid such that all the barrier levels fall on the grid. One technique, described in \citep[p. 171]{TaRa00} consists in using a cubic spline interpolation to map a uniform grid to a smoothly deformed grid. A robust implementation is not so simple, and grid points may end up very close to each other, potentially deteriorating the accuracy of the scheme. Another use case where the smooth deformation may become problematic is the valuation of a portfolio of financial derivative contracts on the same finite difference grid.
 
An alternative technique is to keep the grid simple (for example uniform) and use \emph{ghost points}\footnote{also known as fictitious points} to enforce the boundary condition(s) \citep[p. 1209–1210]{wilmott2013paul}. \citet{healy2022inserting} finds it to be as accurate as the deformed grid approach for standard barrier options contracts, through the use of quadratic interpolation instead of the linear interpolation presented in \citep{wilmott2013paul}.

Here, we take a look at the implications of the ghost point technique in terms of stability for the explicit Euler scheme. In particular, we show that the effective boundary condition is not a Dirichlet condition anymore, and that it significantly reduces the region of stability.

\section{Ghost point for the one-touch option}
We follow \citet{wilmott2013paul} and discretize the problem using a uniform grid between $S_{\min} = 0$ and $S_{\max} = S(0)e^{\left(r(T)-q(T)-\frac{1}{2}\sigma^2(T)\right)T+4\sigma(T)\sqrt{T}}$. While this choice of $S_{\max}$ is not well adapted for a simple one-touch contract, it is relevant for more exotic contracts, where multiple payoffs will be solved on the same grid, such as a vanilla option in the case of a knock-in.
Let $u$ be the index such that $S_{u-1} < L^+ \leq S_u$, with $S_u = u \delta S$, we have
\begin{equation}
\frac{	V^{k+1}_i-V^k_i}{\delta t} = \frac{\sigma_k^2 S_i^2}{2\delta S^2}\left(V_{i+1}^k -2V_i^k+V_{i-1}^k\right) + (r_k-q_k)S_i\frac{V_{i+1}^k-V_{i-1}^k}{2\delta S} - r_k V_i^k\,,
\end{equation}
for $i=1,...,u-1$ where $V^k_i = V(T- k\delta t,S_i)$. This may be rewritten as a tridiagonal system
\begin{equation}
		V^{k+1}_i =   A_{i,i-1}^k V_{i-1}^k+ (1+A_{i,i}^k) V^k_i  +  A_{i,i+1}^k V_{i+1}^k \,.
\end{equation}
with \begin{align*}
	A_{i,i-1}^k=\delta t \left(\frac{\sigma_k^2 S_i^2}{2\delta S^2}-\frac{(r_k-q_k)S_i}{2\delta S} \right)\,,\quad
	A_{i,i}^k =  \delta t \left(-r_k - \frac{\sigma_k^2 S_i^2}{\delta S^2}\right)\,,\quad
	A_{i,i+1}^k = \delta t \left(\frac{\sigma_k^2 S_i^2}{2\delta S^2}+\frac{(r_k-q_k)S_i}{2\delta S} \right) \,.
\end{align*}

The boundary condition at $S_0=0$ reads
\begin{equation*}
\frac{	V^{k+1}_0-V^k_0}{\delta t}= (r_k-q_k)S_0\frac{V_{1}^k-V_{0}^k}{\delta S} - r_kV_0^k\,,
\end{equation*}
where we used a first difference approximation for the first derivative, and assumed the payoff to be linear (or equivalently the diffusion to be zero) at $S=S_0$.

The fictitious value at $S_u$ reads \begin{equation}V_u^{k+1} = \frac{\delta S}{L^+ - S_{u-1}} \times 1 - \frac{S_u-L^+}{L^+ - S_{u-1}}V_{u-1}^{k+1}\, \label{eqn:ghost}\end{equation} using a linear interpolation between $S_{u-1}$ and $S_u$. For $i > u$, we have $V_u^{k+1} = 1$.

The value at $S_u$ is fictitious and may be entirely bypassed by computing directly $V_{u-1}^{k+1}$ through
\begin{align}
	V_{u-1}^{k+1} &= A_{u-1,u-2}^k V_{u-2}^k + (1+A_{u-1,u-1}^k) V^{k}_{u-1} + A_{u-1,u}^k V_{u}^k \nonumber \\
	&=  A_{u-1,u-2}^k V_{u-2}^k + (1+A_{u-1,u-1}^k) V^{k}_{u-1} + A_{u-1,u}^k \left[\frac{\delta S}{L^+ - S_{u-1}} \times 1 -  \frac{S_u-L^+}{L^+ - S_{u-1}}V_{u-1}^{k}\right] \nonumber\\
	&=  A_{u-1,u-2}^k V_{u-2}^k + \left[1+A_{u-1,u-1}^k  - A_{u-1,u}^k\frac{S_u-L^+}{L^+ - S_{u-1}}\right] V^{k}_{u-1}   + A_{u-1,u}^k\frac{\delta S}{L^+ - S_{u-1}} \times 1
\end{align}
In contrast to Equation \ref{eqn:ghost}, this last equation expresses $V_{u-1}^{k+1}$ in terms of $V_{u-2}^k$ and $V_{u-1}^k$ only. It may thus be used to infer the stability condition of the explicit Euler scheme.

Let $\tilde{A}$ be such that $\tilde{A}_{i,j}=A_{i,j}^k$ for $i\in\{0,...,M\}\setminus \{u-1,u\}$ and  $j=0,...,M$. For $i=u-1$, we let $\tilde{A}_{u-1,u-2} = A_{u-1,u-2}^k$ and $\tilde{A}_{u-1,u-1} = A_{u-1,u-1}^k  - A_{u-1,u}^k\frac{S_u-L^+}{L^+ - S_{u-1}}$, $\tilde{A}_{u-1,u}=0$. For $i=u$, we let $\tilde{A}_{u,u}= 1$ and $\tilde{A}_{u,u-1}=\tilde{A}_{u,u+1}=0$.
To guarantee convergence of the explicit scheme, we require $\lVert I+ \tilde{A} \rVert_\infty \leq 1$ . 

We assume the non-diagonal elements to be non-negative\footnote{A necessary but not sufficient condition for the matrix $-A$ to be an M-matrix.}, that is $|(r_k-q_k)|\leq \frac{\sigma^2 S_i}{\delta S}$ for $i=1,...,n$ which is always possible to enforce by reducing $\delta S$ accordingly.
For $i=1,...,u-2$ the norm bound translates to:
		 \begin{align}\begin{cases}1 - \delta t r_k \leq 1\,,&\quad \textmd{ if } 1+A_{i,i} \geq 0\,,\\
		-1 \leq -1 + \delta t \left(r_k + \frac{2\sigma^2 S_i^2}{\delta S^2}\right) \leq 1\,,&\quad \textmd{ if } 1+A_{i,i} < 0\,.\end{cases}\label{eqn:explicit_step_stab_condition2} \end{align}
We will further assume that $r_k \geq 0$. The condition reduces to
	 \begin{align}	\delta t \left(r_k + \frac{\sigma_k^2 S_i^2}{\delta S^2}\right) \leq 1\,.\label{eqn:explicit_step_stab_condition}
	 	 \end{align}
When the barrier level is on the grid, we have $S_u=L^+$ and the  condition at $i=u-1$ also verifies Equation \ref{eqn:explicit_step_stab_condition}.
When the barrier level is not on the grid, the ghost point imposes the additional condition
\begin{align*}
1+A_{u-1,u-1}^k  - A_{u-1,u}^k\frac{S_u-L^+}{L^+ - S_{u-1}} \geq 0\quad \textmd{ and }  1+A_{u-1,u-1}^k  - A_{u-1,u}^k\frac{S_u-L^+}{L^+ - S_{u-1}} + A_{u-1,u-2} \leq 1\,,\quad \textmd{ or }\\
1+A_{u-1,u-1}^k  - A_{u-1,u}^k\frac{S_u-L^+}{L^+ - S_{u-1}} < 0\quad \textmd{ and  }  -1\leq -1-A_{u-1,u-1}^k  + A_{u-1,u}^k\frac{S_u-L^+}{L^+ - S_{u-1}}+ A_{u-1,u-2} \leq 1
\end{align*}
equivalently
\begin{align*}\delta t \left(\frac{\sigma_k^2 S_{u-1}^2}{2\delta S^2}-\frac{(r_k-q_k)S_{u-1}}{2\delta S} \right) \leq	\delta t \left(r_{k} + \frac{\sigma_k^2 S_{u-1}^2}{\delta S^2} +  \left(\frac{\sigma_k^2 S_{u-1}^2}{2\delta S^2}+\frac{(r_k-q_k)S_{u-1}}{2\delta S} \right)\frac{S_u-L^+}{L^+ - S_{u-1}} \right) \leq 1\quad\,,	\textmd{  or } \\
 1 < \delta t \left(r_{k} + \frac{\sigma_k^2 S_{u-1}^2}{\delta S^2} + \left(\frac{\sigma_k^2 S_{u-1}^2}{2\delta S^2}+\frac{(r_k-q_k)S_{u-1}}{2\delta S} \right)\frac{S_u-L^+}{L^+ - S_{u-1}} \right) \leq 2 -\delta t \left(\frac{\sigma_k^2 S_{u-1}^2}{2\delta S^2}-\frac{(r_k-q_k)S_{u-1}}{2\delta S} \right) \,.
 \end{align*}

 
In particular, when $r=q=0$, we obtain \begin{equation}\delta t  \leq \frac{4\delta S^2}{ \sigma^2_k S_{u-1}^2\left(3+\frac{S_u-L^+}{(L^+ - S_{u-1})}\right)}\,.\label{eqn:explicit_step_rates0}\end{equation} The condition becomes increasingly stringeant as the grid point $S_{u-1}$ moves towards $L^+$. Let $\epsilon = L^+-S_{u-1}$, we have \begin{align*}\frac{4\delta S^2}{ \sigma^2_k S_{u-1}^2\left(3+\frac{\delta S-\epsilon}{\epsilon}\right)} = \frac{4\delta S}{ \sigma^2_k S_{u-1}^2}\epsilon+\mathcal{O}(\epsilon^2)\,.\end{align*}

 For any value of $\delta t$, there exists a threshold such that the explicit scheme is unstable as $S_{u-1}$ moves towards $L^+$. Inserting a point in the grid at the barrier level would however be worse since Equation \ref{eqn:explicit_step_stab_condition} leads to a time-step restriction in $\epsilon^2$ instead of $\epsilon$ for the ghost point.
 
 When the barrier falls in the middle of two grid points, the timestep restriction from Equation \ref{eqn:explicit_step_rates0} becomes the same as the one from Equation \ref{eqn:explicit_step_stab_condition}.


\section{Numerical results} 
We consider a one touch option of maturity $T=1$ year and barrier $L^+=7581.36$, on an asset of spot price $S=6317.80$ with rates $r=q=0\%$  and constant volatility $\sigma=20\%$. The reference price  obtained by an analytical formula is 32.9620\%.
\begin{figure}[h]
	\centering{\includegraphics[width=\textwidth]{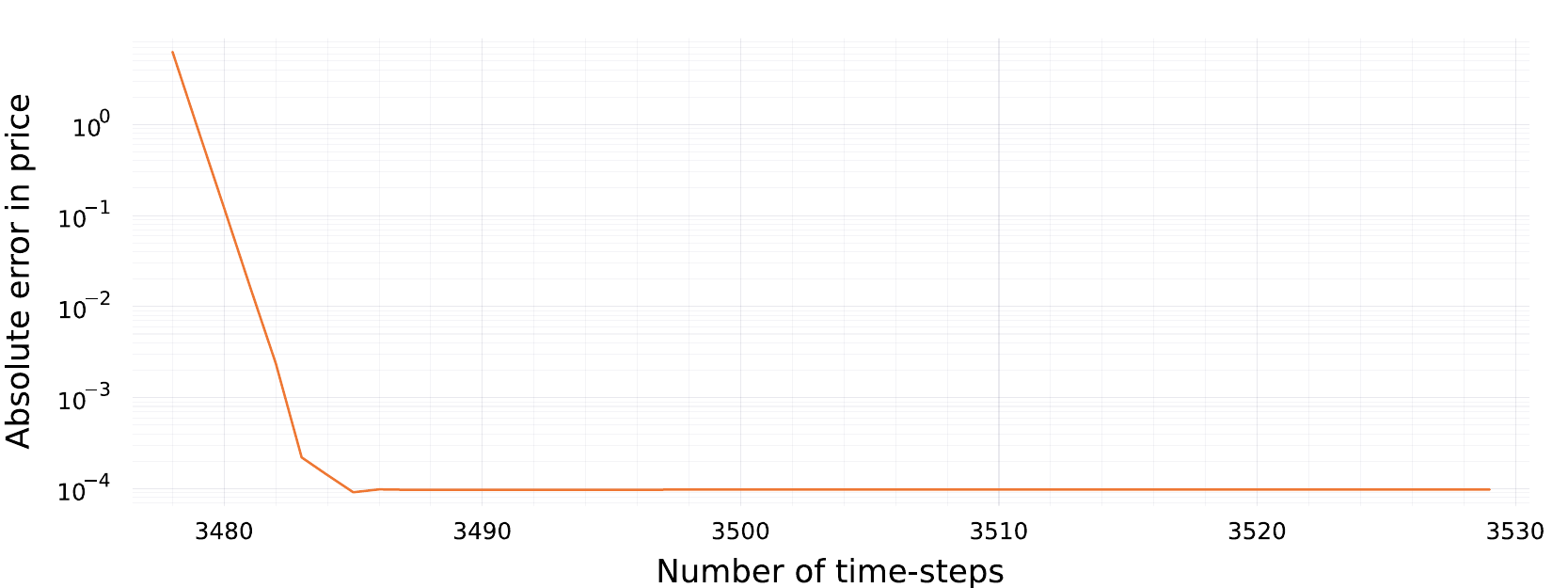}}
	\caption{Error in one-touch option price at $S=S(0)$ as a function of the number of time-steps for a grid composed of $M=100$ space-steps.}
	\label{fig:explicit_rates0_error}
\end{figure}
 We use $M=100$ space-steps with $S_{\min}=0$ and $S_{\max}= 13782.11$ and plot the absolute error in price as a function of the number of time-steps in Figure \ref{fig:explicit_rates0_error}. 
 
\begin{figure}[h]
	\centering{\includegraphics[width=\textwidth]{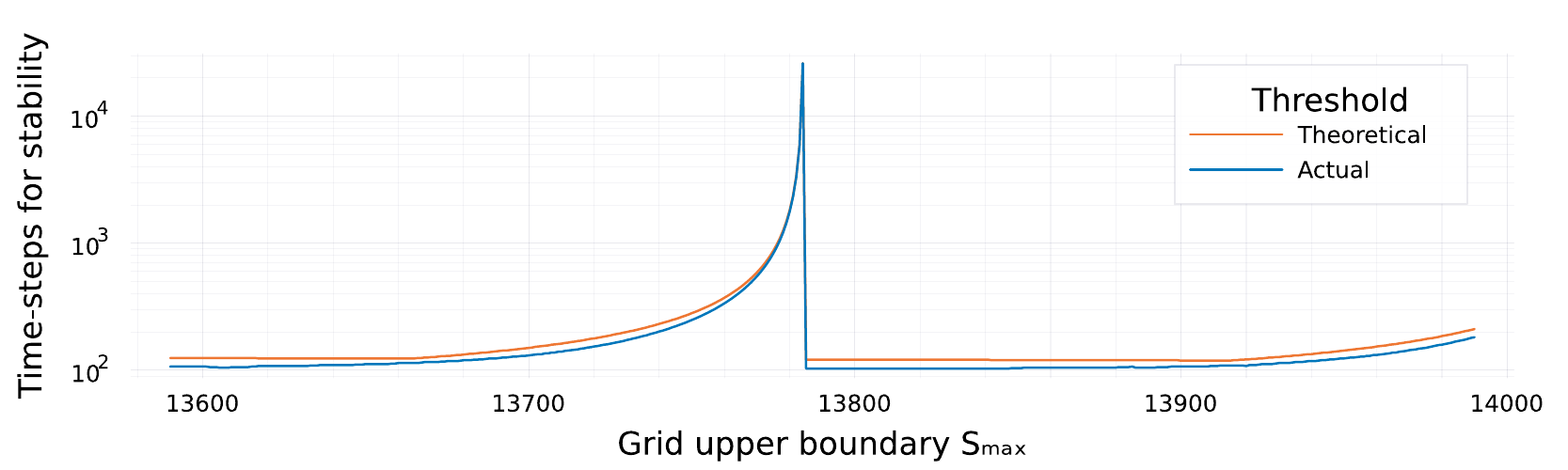}}
	\caption{Theoretical and actual thresholds on the number of time-steps (in log scale) where the one-touch option solution starts diverging for a uniform grid composed of $M=100$ space-steps, varying the grid boundary $S_{\max}$.}
	\label{fig:thresholds}
\end{figure}

In such a setup, the barrier is much closer to $S_{u-1}\approx 7580.16$ than to $S_u\approx 7717.97$.
The threshold from Equation \ref{eqn:explicit_step_rates0} is 3529 while the "standard" explicit timestep threshold from Equation \ref{eqn:explicit_step_stab_condition} is 401, based on $S_{\max}$ and 121, based on $L^+$ (which, on this example is enough, as $V(S)=0$ for $S > L^+$). Stability breaks effectively below around 3485 timesteps. Table \ref{tbl:thresholds} and Figure \ref{fig:thresholds} compare the actual and theoretical thresholds where the numerical solution starts diverging, varying the upper boundary of the grid $S_{\max}$.

\begin{table}
	\centering{
		\caption{Theoretical and actual thresholds on the number of time-steps where the one-touch option solution starts diverging for a uniform grid composed of $M=100$ space-steps, varying the grid boundary $S_{\max}$. \label{tbl:thresholds}.}
	\begin{tabular}{rrrrr}\toprule
$S_{\max}$ & $L^+ - S_{u-1}$ & $(L^+ - S_{u-1}) / \Delta S$ & Theoretical Threshold & Actual Threshold\\\midrule
13662.0 &  67.26 & 0.492 & 122 & 115 \\
13702.0  & 45.26 & 0.330 & 153 & 134 \\
13760.0 & 13.36 & 0.097 & 373 & 338 \\
13772.0 &  6.76 & 0.049 & 677  & 638 \\
13778.0 & 3.46 &  0.021& 1266 & 1223 \\
13782.0 & 1.26 & 0.009&  3370 & 3322\\
13784.0 & 0.16 & 0.001 & 26121 &26071 \\\bottomrule
	\end{tabular}}
\end{table}

\begin{figure}[h]
	\centering{
		\subfloat[][Uniform with ghost point.]{\includegraphics[width=0.5\textwidth]{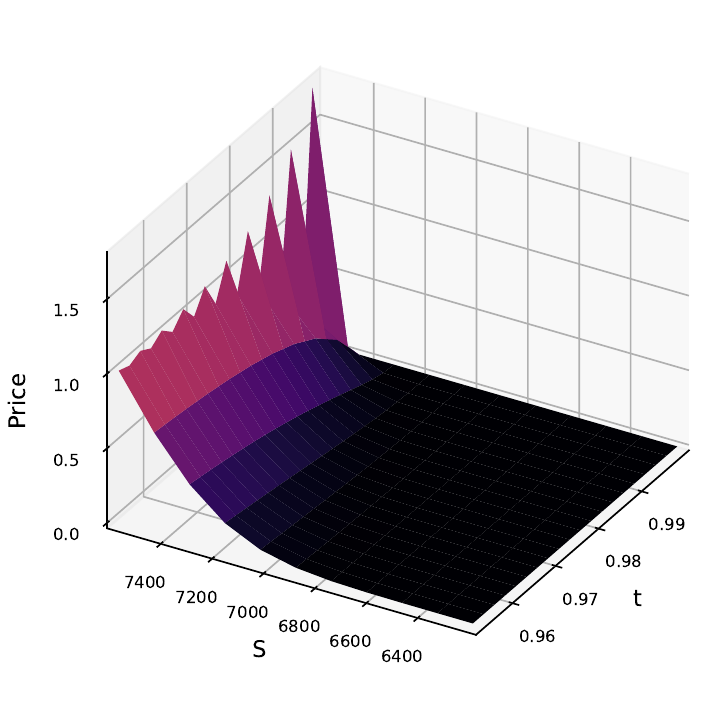}}
	\subfloat[][Streched, barrier on grid.]{\includegraphics[width=0.5\textwidth]{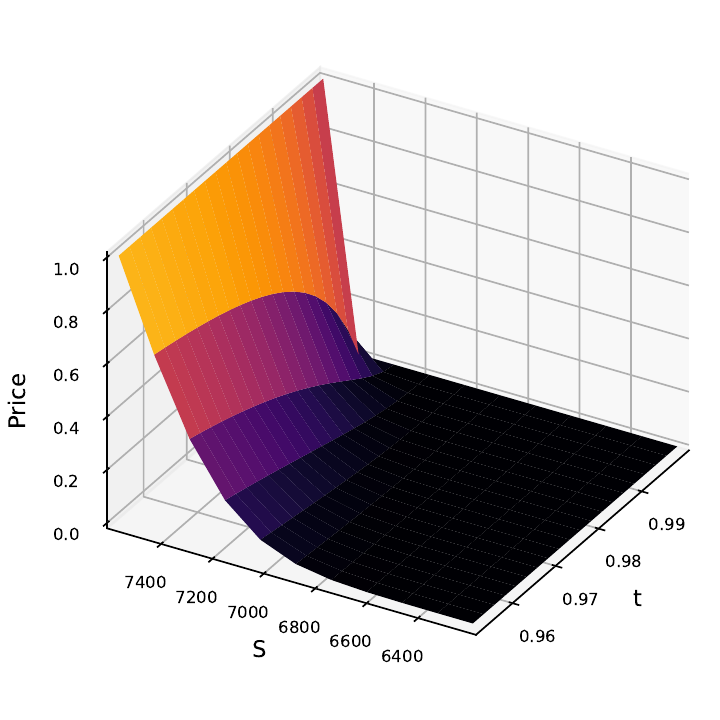}}
	}
	\caption{Prices of a one-touch option obtained by the Crank-Nicolson scheme on the finite difference grid, for $M=100$ space-steps and $N=400$ time-steps, close to the expiry and the barrier level.}
	\label{fig:cn}
\end{figure}

\begin{figure}[h]
	\centering{\includegraphics[width=0.5\textwidth]{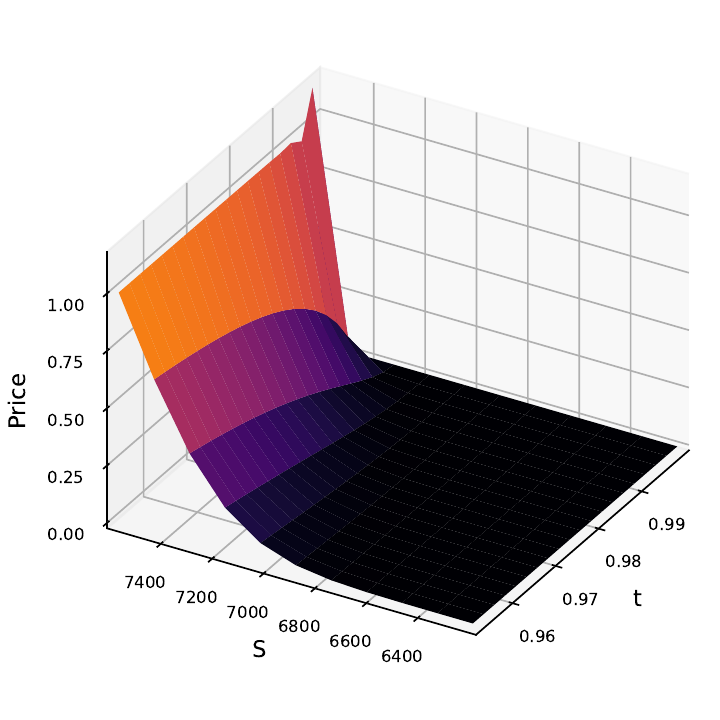}}
	\caption{Prices of a one-touch option obtained by the TR-BDF2 scheme on the finite difference grid, for $M=100$ space-steps and $N=400$ time-steps, close to the expiry and the barrier level on a uniform grid, with ghost point.}
	\label{fig:trbdf2}
\end{figure}

The ghost point boundary condition results in one large eigenvalue for the matrix $I+A$. It may thus create spurious oscillations for some implicit schemes. In Figures \ref{fig:cn}, \ref{fig:trbdf2} we plot the values close to the maturity and close to the barrier level for the Crank-Nicolson scheme and for the TR-BDF2 scheme. 
With 400 time-steps, the Crank-Nicolson scheme preserves the monotonicity of the solution when the barrier is placed on the grid, as expected.
When the barrier is not on the grid, and the ghost point technique is used, a narrow oscillation occurs, limited to the values next to the barrier level for the Crank-Nicolson scheme. The TR-BDF2 scheme damps the oscillations very quickly due to its L-stability properties \citep{lefloch2014tr} and oscillations are not visible on our example (Figure \ref{fig:trbdf2}). 


\section{Conclusion}
The ghost point technique is often not appropriate on the explicit Euler scheme. In the context of a simple one-touch option, the scheme is guaranteed to become unstable as the grid point moves towards the barrier level, regardless of the time-step size.
Inserting a point at the barrier level on the grid is worse, as the divergence will appear for larger time-steps. The correct approach is to place the barrier level on the grid by grid shifting or a smooth grid deformation. In the presence of multiple barrier levels for different time periods, shifting or deformation may be done at a single barrier level, along with a quadratic interpolation at the transioning times where the barrier level changes. Compared to a deformation which places all the barrier levels on the grid, this would avoid the issue of having grid points too close to each other. The shift and interpolate technique is however not applicable for the case of pricing a portfolio of exotic trades on the same grid, the grid must be rebuilt for each trade with a continuous barrier feature.

For L-stable implicit schemes, no oscillation is visible and the ghost point technique is found as accurate as the placement of the barrier level on the grid by a smooth deformation. For the A-stable Crank-Nicolson scheme, a very localized oscillation is present close to the boundary.


\externalbibliography{yes}
\bibliography{explicit_ghost.bib}
\appendixtitles{no}

\end{document}